\documentclass{amsart}
\usepackage{pstricks, pstricks-add, mathptmx, amsfonts, hyperref, color}
\theoremstyle{plain}
\newtheorem{thm}{Theorem}[section]
\newtheorem*{dbp}{$\dbar$-problem for $E$}

\newtheorem{prop}[thm]{Proposition}

\theoremstyle{definition}

\newtheorem{rhp}{Riemann-Hilbert Problem}
\newtheorem*{remark}{Remark}

\newcommand{\dbar}{{\overline{\partial}}}
\newcommand\pmtwo[4]{\left( \begin{array}{cc}#1&#2\\#3&#4\end{array} \right)}

\begin{document}

\title{Long--time Asymptotics for the NLS equation via $\bar{\partial}$ methods}
\author{Momar Dieng}
\address{Dept. of Math, Univ. of Arizona, 520-621-2713, FAX: 520-626-5186}%
\email{momar@math.arizona.edu}%
\author{K. D. T-R McLaughlin}%
\address{Dept. of Math., Univ. of Arizona}%
\email{mcl@math.arizona.edu}

\thanks{The authors were supported in part by NSF
grants DMS-0451495 and DMS-0800979.}

\begin{abstract}
We present a new method for obtaining sharp asymptotics of solutions of the defocussing nonlinear Schr\"odinger (NLS) equation, based on $\dbar$ methods and under essentially minimal regularity assumptions on initial data.
\end{abstract}

\maketitle

\section{Introduction}

The long time behavior of solutions $q(x,t)$ of the defocusing Nonlinear Schr\"{o}dinger (NLS) equation
\begin{eqnarray}
\label{eq:NLSEQ}
\left\{ \begin{array}{c}
i q_{t} + q_{xx} - 2 \left| q \right|^{2} q = 0, \\
q(x,t=0) = q_{0}(x) \rightarrow 0 \mbox{ as } |x| \to \infty \\
\end{array} \right.
\end{eqnarray}
has been studied extensively, under varying degree of smoothness assumptions on the initial data $q_{0}$ \cite{Zakh1, Deif3, Deif4, Deif7, Deif5}.  The asymptotic behavior takes the following form: as $t \to \infty$, one has
\begin{eqnarray}
\label{eq:AsForm}
q(x,t) = t^{-1/2}\alpha(z_0)e^{ix^2/(4t) - i\nu(z_0)\log(8t)} + \mathcal{E}\left(x,t\right),
\end{eqnarray}
where
\[
\nu(z) = -\frac1{2\pi}\log(1-|r(z)|^2),\quad |\alpha(z)|^2=\nu(z)/2,
\]
and
\[
\hspace{-0.3in}\arg\alpha(z)=\frac1\pi\int_{-\infty}^{z}\log(z-s)d(\log(1-|r(s)|^2))+\frac\pi4 + \arg\Gamma(i\nu(z))-\arg r(z).
\]
Here $z_{0} = - x / (4t)$, $\Gamma$ is the gamma function and the function $r$ is the so-called reflection coefficient associated to the initial data $q_{0}$, as described later in this section.

Estimates on the size of the error term $\mathcal{E}\left(x, t \right)$ depend on the smoothness assumptions on $q_{0}$. The above asymptotic form was first obtained in \cite{Zakh1}.  The nonlinear steepest descent method \cite{Deif6} was brought to bear on this problem in \cite{Deif3} (see \cite{Deif4} for a pedagogic description), where the authors assumed the initial data possessed high orders of smoothness and decay, and proved that $\mathcal{E}\left(x,t\right)$ satisfied
\begin{eqnarray}
 \mathcal{E}\left(x,t\right) = \mathcal{O} \left( \frac{\log{t}}{t} \right)
\end{eqnarray}
uniformly for all $x \in \mathbb{R}$.

Early in this millennium, Deift and Zhou developed some new tools for the analysis of Riemann--Hilbert problems, originally aimed at studying the long time behavior of perturbations of the NLS equation \cite{Deif8}.  Their methods allowed them to establish long time asymptotics for the NLS equation with essentially minimal assumptions on the initial data \cite{Deif5}.  Indeed, they showed that if the initial data is in the Sobolev space \[H^{1,1}=\{f\in L^2(\mathbb R):\ xf,f'\in L^2(\mathbb R)\},\]
then the (unique, weak) solution satisfies (\ref{eq:AsForm}) with an error term $\mathcal{E}\left(x,t \right)$  that satisfies, for any fixed $0 < \kappa < 1 / 4$,
\begin{eqnarray}
 \mathcal{E}\left(x,t \right)  = \mathcal{O}\left( t^{- \left( \frac{1}{2} + \kappa \right) } \right) \ .
\end{eqnarray}

Recently, McLaughlin and Miller \cite{McLa1,Mcla2} have developed a method for the asymptotic analysis of Riemann--Hilbert problems based on the analysis of $\dbar$ problems, rather than the asymptotic analysis of singular integrals on contours.  In this paper we adapt and extend this method to the Riemann--Hilbert problem associated to the NLS equation.  The main point of our work is this:  by using the $\dbar$ approach, we avoid all delicate estimates involving Cauchy projection operators in $L^{p}$ spaces (which are central to the work in \cite{Deif5}).  In place of such estimates, we carry out basic estimates of double integrals which involve nothing more than calculus.  Our result is as follows.

\begin{thm}
For initial data $q_0(x)=q(x,t=0)$ in the Sobolev space \[H^{1,1}=\{f\in L^2(\mathbb R):\ xf,f'\in L^2(\mathbb R)\},\] we have as $t\to\infty$:
\[
q(x,t) = t^{-1/2}\alpha(z_0)e^{ix^2/(4t) - i\nu(z_0)\log(8t)} + \mathcal{O}\left({t^{-3/4}}\right)\ .
\]
where
\[
\nu(z) = -\frac1{2\pi}\log(1-|r(z)|^2),\quad |\alpha(z)|^2=\nu(z)/2,
\]
and
\[
\hspace{-0.3in}\arg\alpha(z)=\frac1\pi\int_{-\infty}^{z}\log(z-s)d(\log(1-|r(s)|^2))+\frac\pi4 + \arg\Gamma(i\nu(z))-\arg r(z).
\]
\label{mainresult}\end{thm}
The main features of this result are threefold:
\begin{enumerate}
\item the error term is an improvement (see \cite{Deif5}; in fact our estimate on the error is sharp)
\item the new $\dbar$ method which is used to derive it affords a considerably less technical proof than previous results.
\item the method used to establish this result is readily extended to derive a more detailed asymptotic expansion, beyond the leading term (see the remark at the end of the paper).
\end{enumerate}

\vskip 0.2in
The solution procedure for the nonlinear Schr\"{o}dinger equation may be described as follows:  for a given function $r(z)$ in the Sobolev space $H^{1,1}_{1} = \{ f \in L^{2}(\mathbb{R}): z f, f' \in L^{2}(\mathbb{R}), \sup_{z \in \mathbb{R}} |f(z)| < 1 \}$, consider the following Riemann--Hilbert problem:
\begin{rhp}
\label{rhp:01}
Find $M = M(z) = M(z;x,t)$ a $2 \times 2$ matrix, satisfying the following conditions:
\begin{equation}
\begin{cases}
M \quad\textrm{analytic on $\mathbb C\setminus\Sigma$}\\
M_{+}(z)=M_{-}(z)V_{M}(z)\quad\textrm{for $z\in\Sigma=\mathbb R$ with $V_{M}$ specified below}\\
M= I+O\left(\frac{1}{z}\right) \quad\textrm{as $z\to\infty$}.
\end{cases}
\end{equation}
The jump  matrix $V_{M}$ is defined on $\Sigma$ as follows:
\[
V_{M}=\left(
  \begin{array}{cc}
    1-|r|^{2} & -\bar{r}\,e^{-2it\theta} \\
r\,e^{2it\theta} & 1
  \end{array}
\right), \qquad \theta=2z^{2}-4z_{0}z = 2 z^{2} + x z / t.
\]
\end{rhp}
Next, define
\begin{eqnarray}
M_{1}(x,t) = \lim_{z \to \infty} z \left( M(z) - I \right),
\end{eqnarray}
and then set
\begin{eqnarray}\label{qMrelation}
q(x,t) = 2 i \left( M_{1}(x,t) \right)_{12}.
\end{eqnarray}
The fact of the matter is that $q(x,t)$ solves the nonlinear Schr\"{o}dinger equation (\ref{eq:NLSEQ}).  The connection between the initial data $q_{0}(x)$ and the reflection coefficient $r(z)$ is achieved through the spectral and inverse-spectral theory of the associated Zakharov-Shabat differential operator
\begin{equation}
\mathcal{L}=
\nonumber
i \left( \begin{array}{cc}
1 & 0 \\
0 & -1 \\
\end{array} \right) \frac{d}{dx} + \left( \begin{array}{cc}
0 & - i q \\
i \overline{q} & 0 \\
\end{array} \right),
\end{equation}
as described, for example, in \cite{Deif4}.  It is well known that if $q_{0} \in H^{1,1}$, then $r \in H^{1,1}_{1}$, and more generally this spectral transform, $\mathcal{R}: H^{1,1} \to H^{1,1}_{1}$, $q_{0} \mapsto r$, is a bijection between these two spaces (it turns out that $\mathcal{R}$ is bi-Lipschitz) \cite{Zhou1}. For initial data in $H^{1,1}$, $q(x,t)$ obtained from the Riemann--Hilbert problem described above is the unique weak solution to the NLS equation (see \cite{Deif5}).

\vskip 0.2in
Recent analyses of the long-time behavior of the solution of the NLS initial value problem (\ref{eq:NLSEQ}) have involved the detailed analysis of the behavior of the solution $M$ to the Riemann--Hilbert problem \ref{rhp:01}.  As regularity assumptions on the initial data $q_{0}$ are relaxed, the detailed analysis becomes more involved, technically.  The purpose of this manuscript is to carry out a complete analysis of the long-time asymptotic behavior of $M$ under the assumption that $r \in H^{1,1}_{1}$, as in \cite{Deif4}, but via a $\overline{\partial}$ approach which replaces very technical harmonic analysis involving Cauchy projection operators with very straightforward estimates involving some explicit two-dimensional integrals.

\section{Proof of the result}

We outline the proof of the results. We begin by outlining the solutions to two model RHPs which we shall need in our derivation.

\subsection{A model scalar RHP}

Let $z_{0}\in\mathbb R$, and $\sigma_{3}=\pmtwo{1}{0}{0}{-1}$. We seek a function $\delta(z)$ satisfying:
  \begin{equation}
  \begin{cases}
   \textrm{ $\delta$ is analytic and invertible for $z \in \mathbb{C} \setminus {\mathbb{R}},$}\\
   \textrm{$\delta(z) \to 1$ as $z\to\infty$,}\\
   \delta_{+}(z) = \begin{cases} \delta_{-}(z) \left( 1 - \left| r(z) \right|^{2}\right), &  z < z_{0} \\ \delta_{-}(z) & z> z_{0}\end{cases} .
  \end{cases}
  \label{modrhp1}
  \end{equation}

The unique solution to this RHP is easily checked to be the function
\begin{equation}
\delta(z) = \exp[\gamma(z)]=\exp\left[\frac{1}{2i\pi}\int_{-\infty}^{z_{0}}\frac{\ln\left(1-|r(s)|^{2}\right)}{s-z}\,ds\right], \qquad z\notin\mathbb R
\label{modrhp1sol}.
\end{equation}
We will need estimates on $\delta(z)$ which we take from \cite{Deif5} (see Proposition~(2.12) and formula~(2.41) of that paper).
Suppose $r\in L^\infty({\mathbb R}) \cap L^2({\mathbb R})$ and $\|r\|_{L^\infty} \le \rho <1$. Then
\begin{equation}
\delta(z) \overline{\delta(\bar z)} = 1,
\end{equation}
\begin{equation}
(1-\rho)^{\frac12} \le (1-\rho^2)^{\frac12} \le |\delta(z)|, |\delta^{-1}(z)| \le (1-\rho^2)^{-\frac12} \le (1-\rho)^{-\frac12},\label{deltabound}
\end{equation}
and
\begin{equation}|\delta^{\pm 1}(z)|\le 1 \quad \text{for}\quad
\pm \text{Im } z>0.
\end{equation}
For real $z$,
\begin{equation}
|\delta_+(z) \delta_-(z)| = 1 \quad \text{and, in
particular,}\quad |\delta(z)|  = 1\quad \text{for} \quad
z>z_0,
\end{equation}
\begin{equation}|\delta_+(z)| = |\delta^{-1}_-(z)| = (1 - |r(z)|^2
)^{\frac12},\qquad z<z_0,
\end{equation}
and
\begin{equation}\Delta\equiv \delta_+\delta_- = e^{\frac1{i\pi}
\text{ P.V. } \int^{z_0}_{-\infty}
\frac{\log(1-|r(s)|^2)}{s-z}ds}, \quad \text{where P.V.\ denotes
the principal value.}
\end{equation}
Also $|\Delta| = |\delta_+\delta_-| = 1$, and
\begin{equation}
\|\delta_\pm - 1\|_{L^2(dz)} \le \frac{c\|r\|_{L^2}}{1-\rho}.
\end{equation}
Let $\chi^0(s)$ denote the
characteristic function of the interval $(z_0-1,z_0)$. Then for $z\in\mathbb C\setminus~(-\infty,z_0)$,
\begin{align}
\gamma(z) &= \int^{z_0}_{-\infty} \{\log(1-|r(s)|^2) - \log(1-|r(z_0)|^2) \chi^0(s)(s-z_0+1)\} \frac{ds}{2\pi i(s-z)}\label{gammaexpansion}\\
&\quad + \log(1-|r(z_0)|^2) \int^{z_0}_{z_0-1} \frac{s-z_0+1}{s-z} \frac{ds}{2\pi i}\nonumber\\
&= \beta(z,z_0) + i\nu(z_0)(1-(z-z_0 + 1) \log(z-z_0+1)\nonumber\\
&\quad + (z-z_0) \log(z-z_0) + \log(z-z_0)),\nonumber
\end{align}

where $\nu(z_0) = -\frac1{2\pi} \log(1-|r(z_0)|^2)$, and
\begin{equation}
\beta(z,z_0) = \int^{z_0}_{-\infty} \{\log(1-|r(s)|^2) -
\log(1-|r(z_0)|^2) \chi^0(s) (s-z_0+1)\} \frac{ds}{2\pi i(s-z)}.
\label{beta}\end{equation}
On any ray $L_{\phi} = z_0 + e^{-i\phi} \mathbb R_{+} =\{z =
z_0 + ue^{-i\phi}, u\ge 0\}$ with $0<\phi<\pi$ or $-\pi<\phi<0$ we find
\begin{equation}
\|\beta\|_{H^1(L_{\phi})} \le \frac{\hat{c}\|r\|_{H^{1,0}}}{1-\rho},
\label{betaestimate1}\end{equation}
where the constant $\hat{c}$ is independent of $\phi$ and $z_{0}$ (see Lemma~23.3 in \cite{Beal1}). Then, by standard Sobolev estimates, it follows that there is a constant $c$ independent of $\phi$ and $z_{0}$ such that on $L_{\phi}$
\begin{equation}
\begin{cases}\beta(z,z_0)\textrm{ is continuous up to $z=z_0$}\\
\|\beta(.,z_0)\|_{L^\infty(L_{\phi})} \le c\|r\|_{H^{1,0}}/(1-\rho) \\
|\beta(z,z_0) - \beta(z_0,z_0)| \le \frac{c\|r\|_{H^{1,0}}}{1-\rho} |z-z_0|^{1/2}.
\end{cases}
\label{betaestimate2}\end{equation}
Notice that (\ref{betaestimate2}) also provides us with estimates that are uniformly valid in $\mathbb C^{+}$ and $\mathbb C^{-}$ separately.

\subsection{A model matrix RHP}

We will refer to the RHP described in this section as the Parabolic Cylinder RHP for reason that will soon be obvious. Given $r_{0}$, a complex number satisfying $|r_{0}|<1$, set $\nu = \frac{-1}{2 \pi } \log{\left( 1 - \left| r_{0} \right|^{2}\right)}$, and then seek $P(\xi)$, a $2 \times 2$ matrix
valued function of $\xi$, satisfying:
\begin{equation}
\begin{cases}
\textrm{$P$ is analytic for $\xi \in \mathbb{C} \setminus \mathbb{R}$},\\
P(\xi) = \mathbb{I} + \frac{P_{1}^{\infty}}{\xi} +\mathcal{O}  \left( \frac{1}{\xi^{2}}\right),\qquad \xi\to\infty\\
\textrm{$P_{+}(\xi) = P_{-}(\xi) V_{P}(\xi)$, for $\xi \in \Sigma_{P}$}.
\end{cases}
\label{modrhp2}
\end{equation}
The contour $\Sigma_{P}$ consists of four rays emanating from the
origin, one in each of the four quadrants and described parametrically
as follows:
\begin{eqnarray*}
&& \Sigma_{P}^{I}= \{ \xi: \ \xi = r e^{i \pi / 4},\ r \ge 0 \} , \\
&& \Sigma_{P}^{II}= \{ \xi: \ \xi = r e^{3 i \pi / 4},\ r \ge 0 \}, \\
&& \Sigma_{P}^{III}= \{ \xi: \ \xi = r e^{5 i \pi / 4},\ r \ge 0 \}, \\
&& \Sigma_{P}^{IV}= \{ \xi: \ \xi = r e^{7 i \pi / 4},\ r \ge 0 \}.
\end{eqnarray*}
The jump matrix $V_{P}$ is defined separately on each of the four
rays, as follows:
\begin{eqnarray*}
V_{P}(\xi) := \left\{ \begin{array}{cc}
\pmtwo
{1}{0}{r_{0}\xi^{-2 i \nu}e^{i \xi^{2}/2}}{1}& \ \mbox{ for } \xi \in \Sigma_{P}^{I}\\
\pmtwo
{1}{\frac{-\overline{r_{0}}}{1 - \left |r_{0} \right|^{2}}\xi^{2 i \nu}e^{-i \xi^{2}/2}}{0}{1}& \ \mbox{ for } \xi \in \Sigma_{P}^{II}\\
\pmtwo
{1}{0}{\frac{r_{0}}{1 - \left |r_{0} \right|^{2}}
\xi^{-2 i \nu}e^{i \xi^{2}/2}}{1}& \ \mbox{ for } \xi \in \Sigma_{P}^{III}\\
\pmtwo
{1}{-\overline{r_{0}}\xi^{2 i \nu}e^{-i \xi^{2}/2}}{0}{1}& \ \mbox{ for } \xi \in \Sigma_{P}^{IV}
\end{array}
\right\}
\end{eqnarray*}

\bigskip
\begin{figure}[ht]
\begin{center}
\begin{pspicture}(-3,-3)(3,3)
\rput(0,0){\rnode{0}{}}
\rput(3,3){\rnode{1}{}}
\rput(-3,3){\rnode{2}{}}
\rput(-3,-3){\rnode{3}{}}
\rput(3,-3){\rnode{4}{}}
\ncline[ArrowInside=->,arrowscale=2]{0}{1}
\ncline[ArrowInside=->,arrowscale=2]{3}{0}
\ncline[ArrowInside=->,arrowscale=2]{2}{0}
\ncline[ArrowInside=->,arrowscale=2]{0}{4}
\psset{linestyle=dashed}
\rput(0,0){\psline(-4,0)(4,0)}
\rput(1.2,0.5){$\Omega_{1}$}
\rput(0,0.8){$\Omega_{2}$}
\rput(-1.2,0.5){$\Omega_{3}$}
\rput(-1.2,-0.5){$\Omega_{4}$}
\rput(0,-0.8){$\Omega_{5}$}
\rput(1.2,-0.5){$\Omega_{6}$}
\rput(4.2,1.7){
$V_{P}^{I}=\pmtwo{1}{0}{r_{0}\xi^{-2 i \nu}e^{i \xi^{2}/2}}{1}$
}
\rput(4.2,-1.7){
$V_{P}^{IV}=\pmtwo{1}{-\overline{r_{0}}\xi^{2 i \nu}e^{-i \xi^{2}/2}}{0}{1}$
}
\rput(-4.5,-1.7){
$V_{P}^{III}=\pmtwo{1}{0}{\frac{r_{0}}{1 - \left |r_{0} \right|^{2}}\xi^{-2 i \nu}e^{i \xi^{2}/2}}{1}$
}
\rput(-4.5,1.7){
$V_{P}^{II}=\pmtwo{1}{\frac{-\overline{r_{0}}}{1 - \left |r_{0} \right|^{2}}\xi^{2 i \nu}e^{-i \xi^{2}/2}}{0}{1}$
}
\end{pspicture}
\caption{The countour $\Sigma_{P}$ and the regions $\Omega_{i}, i=1,\ldots,6$.}
\label{sigmaPC}
\end{center}
\end{figure}
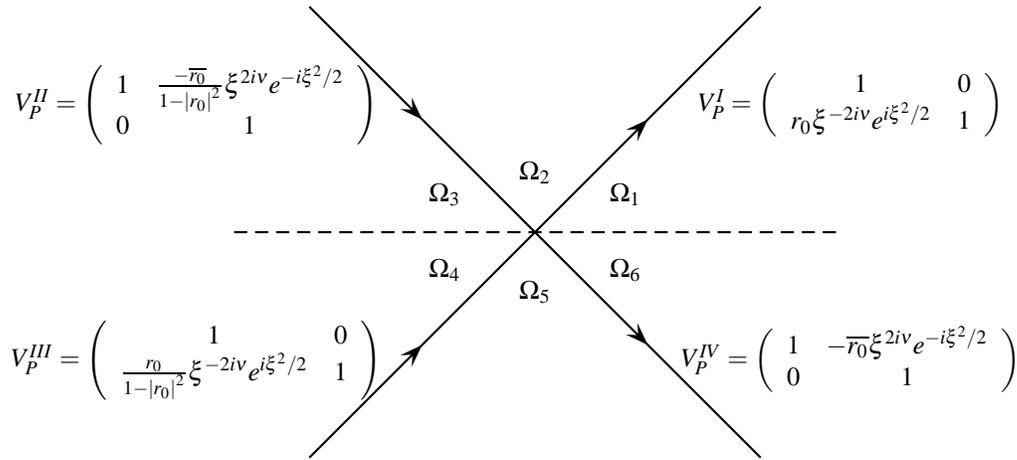
The solution to this Riemann-Hilbert problem is the explicit piecewise analytic matrix valued function $P(\xi)$, defined below \cite{Deif3, Deif4}.  We first define two auxiliary matrix valued functions $\Psi^{+}(\xi)$ and $\psi^{-}(\xi)$, defined in $\mathbb{C}_{+}$ and $\mathbb{C}_{-}$ respectively, as follows.
For $\xi \in \mathbb{C}_{+}$
{\scriptsize
\[\Psi^{+}(\xi) = \pmtwo{e^{- 3 \pi \nu / 4} D_{i \nu}(e^{-3 i \pi /4} \xi)}
{e^{\pi \nu / 4}\left( \beta_{21}\right)^{-1}  \left[ \partial_{\xi} \left( D_{-i \nu}(e^{- i \pi / 4}\xi)\right) - \frac{i \xi}{2} D_{-i \nu}(e^{-  i \pi / 4}\xi) \right]}
{e^{-3 \pi \nu / 4}\left( \beta_{12}\right)^{-1}  \left[ \partial_{\xi} \left( D_{i \nu}(e^{-3 i \pi / 4}\xi)\right) + \frac{i \xi}{2} D_{i \nu}(e^{- 3 i \pi / 4}\xi) \right]}
{e^{\pi \nu / 4} D_{-i \nu}(e^{- i \pi / 4}\xi)}. \]}
For $\xi \in \mathbb{C}_{-}$
{\scriptsize
\[\Psi^{-}(\xi) =
\pmtwo{e^{ \pi \nu / 4} D_{i \nu}(e^{ i \pi /4} \xi)}
{e^{-3 \pi \nu / 4}\left( \beta_{21}\right)^{-1}  \left[ \partial_{\xi} \left( D_{-i \nu}(e^{3i \pi / 4}\xi)\right) - \frac{i \xi}{2} D_{-i \nu}(e^{3  i \pi / 4}\xi) \right]}
{e^{\pi \nu / 4}\left( \beta_{12}\right)^{-1}  \left[ \partial_{\xi} \left( D_{i \nu}(e^{ i \pi / 4}\xi)\right) + \frac{i \xi}{2} D_{i \nu}(e^{ i \pi / 4}\xi) \right]}
{e^{-3 \pi \nu / 4} D_{-i \nu}(e^{3 i \pi / 4}\xi)}.\]
}

\vspace{0.2in}
For a fixed complex number $a$, the function $D_{a}(\zeta)$ is a special function, defined to be the unique function satisfying the parabolic cylinder equation
\[\frac{d^{2}}{d\zeta^{2}} D_{a}(\zeta) + \left( \frac{1}{2} - \frac{\zeta^{2}}{4} + a \right) D_{a}(\zeta)= 0,\]
with boundary condition
\[
D_{a}(\zeta)
=\begin{cases}
 \zeta^{a} e^{- \zeta^{2}/4} \left( 1 + \mathcal{O} \left( \zeta^{-2}\right) \right), \\
 \qquad\qquad\textrm{for $\zeta \rightarrow \infty, \ \left| \mbox{arg} \zeta \right| < \frac{3 \pi }{4}$},\\
 \zeta^{a} e^{- \zeta^{2}/4} \left( 1 + \mathcal{O} \left( \zeta^{-2}\right) \right)-(2\pi)^{1/2}\left(\Gamma(-a)\right)^{-1}e^{a\pi\,i}\zeta^{-a-1} e^{\zeta^{2}/4} \left( 1 + \mathcal{O} \left( \zeta^{-2}\right) \right), \\
 \qquad\qquad\textrm{for $\zeta \rightarrow \infty, \frac{\pi}{4}< \arg \zeta < \frac{5 \pi }{4}$},\\
 \zeta^{a} e^{- \zeta^{2}/4} \left( 1 + \mathcal{O} \left( \zeta^{-2}\right) \right)-(2\pi)^{1/2}\left(\Gamma(-a)\right)^{-1}e^{-a\pi\,i}\zeta^{-a-1} e^{\zeta^{2}/4} \left( 1 + \mathcal{O} \left( \zeta^{-2}\right) \right), \\
 \qquad\qquad\textrm{for $\zeta \rightarrow \infty, -\frac{5\pi}{4}< \arg \zeta < -\frac{\pi }{4}$}.
\end{cases}
\]
The parameters $\beta_{21}$ and $\beta_{12}$ are defined as follows:
\begin{eqnarray*}
&&
\beta_{12} = \frac{(2 \pi)^{1/2} e^{i \pi / 4 }e^{- \pi \nu /2 }}
{r_{0} \Gamma(-i\nu)}, \\
&&
\beta_{21}=\frac{(2 \pi)^{1/2} e^{i \pi / 4 }e^{- \pi \nu /2 }}
{r_{0}\nu \Gamma(-i\nu)}.
\end{eqnarray*}

The matrix $P$ is defined in six non-overlapping regions:
\begin{eqnarray*}
\Omega_{1} = \left\{ \xi:  \mbox{arg}(\xi) \in \left(0, \frac{\pi}{4}\right) \right\},
\ \ \ \ &&
\Omega_{2} = \left\{ \xi:  \mbox{arg}(\xi) \in \left(\frac{\pi}{4}, \frac{3\pi}{4}\right) \right\},\\
\Omega_{3} = \left\{ \xi:  \mbox{arg}(\xi) \in \left(\frac{3\pi}{4}, \pi \right) \right\},
\ \ \ \  &&
\Omega_{4} = \left\{ \xi:  \mbox{arg}(\xi) \in \left(\pi, \frac{5\pi}{4}\right) \right\}, \\
\Omega_{5} = \left\{ \xi:  \mbox{arg}(\xi) \in \left(\frac{5\pi}{4}, \frac{7\pi}{4}\right) \right\},
\ \ \ \ &&
\Omega_{6} = \left\{ \xi:  \mbox{arg}(\xi) \in \left(\frac{7\pi}{4}, 2\pi \right) \right\}.
\end{eqnarray*}

Now $P(\xi)$ is defined in each of these regions as follows (see Fig~\ref{sigmaPC}.):

\begin{eqnarray*}
\mbox{For }\xi \in \Omega_{1}: &&P(\xi) =\
\Psi^{+} (\xi)\ \ e^{ i \xi^{2} \sigma_{3}/4} \ \xi^{- i \nu \sigma_{3}}  \
\pmtwo{1}{0}{-r_{0}}{1}
\\
\mbox{For }
\xi \in \Omega_{2}: &&P(\xi) = \ \Psi^{+}(\xi)\  e^{ i \xi^{2} \sigma_{3}/4} \ \xi^{- i \nu \sigma_{3}}
\\
\mbox{For } \xi \in \Omega_{3}: && P(\xi) = \ \Psi^{+}(\xi)\  e^{ i \xi^{2} \sigma_{3}/4} \ \xi^{- i \nu \sigma_{3}}  \
\pmtwo{1}{\frac{\overline{r_{0}}}{1 - \left| r_{0}\right|^{2}}}{0}{1}
\\
\mbox{For } \xi \in \Omega_{4}: && P(\xi) = \ \Psi^{-}(\xi)\ e^{ i \xi^{2} \sigma_{3}/4} \ \xi^{- i \nu \sigma_{3}} \
\pmtwo{1}{\frac{r_{0}}{1 - \left| r_{0}\right|^{2}}}{0}{1}
\\
\mbox{For } \xi \in \Omega_{5}: && P(\xi) = \ \Psi^{-}(\xi) \ e^{ i \xi^{2} \sigma_{3}/4} \ \xi^{- i \nu \sigma_{3}}
\\
\mbox{For } \xi \in \Omega_{6}: &&  P(\xi) = \ \Psi^{-}(\xi) \ e^{ i \xi^{2} \sigma_{3}/4} \ \xi^{- i \nu \sigma_{3}} \
\pmtwo{1}{-\overline{r_{0}}}{0}{1}
\label{PCsol}.\end{eqnarray*}

The fact is that $P$ defined above solves RHP~(\ref{modrhp2}). Moreover in the notation of RHP~(\ref{modrhp2}), we have (see \cite{Deif4})
\begin{equation}P_{1}^{\infty}=\pmtwo{0}{\frac{-i(2\pi)^{1/2}e^{i\pi/4}e^{-\pi\nu/2}}{r_{0}\Gamma(-i\nu)}}{\frac{i\nu r_{0}\Gamma(-i\nu)}{(2\pi)^{1/2}e^{i\pi/4}e^{-\pi\nu/2}}}{0}.\label{p1infinty}\end{equation}
Note also that in the setup of RHP~(\ref{modrhp2}), the origin is the reference point from which the rays emanate. However in the following sections we will need to use the rescaling $\xi\to\sqrt{8t}(z-z_{0}) $ for $z_{0}\in \mathbb R$.

\subsection{The RHP for NLS}

Let $\Sigma=\mathbb R$. We seek a matrix $M$ satisfying the following conditions:
\begin{equation}
\begin{cases}
M \quad\textrm{analytic on $\mathbb C\setminus\Sigma$}\\
M_{+}(z)=M_{-}(z)V_{M}(z)\quad\textrm{for $z\in\Sigma=\mathbb R$ with $V_{M}$ specified below}\\
M= I+O\left(\frac{1}{z}\right) \quad\textrm{as $z\to\infty$}.
\end{cases}\label{nlsrhp}
\end{equation}
The jump  matrix $V_{M}$ defined on $\Sigma$ as follows:
\[
V_{M}=\left(
  \begin{array}{cc}
    1-|r|^{2} & -\bar{r}\,e^{-2it\theta} \\
r\,e^{2it\theta} & 1
  \end{array}
\right), \qquad \theta=2z^{2}-4z_{0}z.
\]
\subsubsection{Step 1: jump matrix factorization}

To the left of $z_{0}$ we use the factorization
\[
V_{M}=\left(
  \begin{array}{cc}
    1 & 0 \\
\frac{r}{1-|r|^{2}}\,e^{2it\theta} & 1
  \end{array}
\right)
\cdot
\left(
  \begin{array}{cc}
    1-|r|^{2} & 0 \\
0 & \frac{1}{1-|r|^{2}}
  \end{array}
\right)
\cdot
\left(
  \begin{array}{cc}
    1 & \frac{-\bar{r}}{1-|r|^{2}}\,e^{-2it\theta} \\
0 & 1
  \end{array}
\right) = U_{L}U_{0}U_{R}
\]
To the right of $z_{0}$ we use the factorization
\[
V_{M}=\left(
  \begin{array}{cc}
    1 & -\bar{r}\,e^{-2it\theta} \\
0 & 1
  \end{array}
\right)
\cdot
\left(
  \begin{array}{cc}
    1 & 0 \\
r\,e^{2it\theta} & 1
  \end{array}
\right) = W_{L}W_{R}
\]
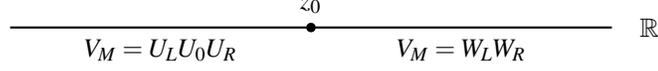
\begin{figure}[ht]
\begin{center}
\begin{pspicture}(-1,-1)(1,1)
\rput(0,0){\psline(0,0)(-4,0)}
\rput(0,0){\psline{*-}(0,0)(4,0)}
\rput(0,0.3){$z_{0}$}
\rput(-2,-0.3){$V_{M}=U_{L}U_{0}U_{R}$}
\rput(2,-0.3){$V_{M}=W_{L}W_{R}$}
\rput(4.5,0){$\mathbb R$}
\end{pspicture}
\caption{The original contour $\mathbb R$ together with the jump matrix in its factorized forms.}
\end{center}
\end{figure}
\subsubsection{Step 2: Extension of $r$ and contour deformation}

We would like to deform the contour by opening sectors, so we need extensions off $\mathbb R$ of the off-diagonal functions $r$, $\bar{r}$,$\frac{r}{1-|r|^{2}}$ and $\frac{-\bar{r}}{1-|r|^{2}}$. We choose these extension in a way that that will make precise later. Let
\begin{enumerate}
\item $R_{1}$ be the extension of $r$ in the sector $\Omega_{1}:\{z:0<\arg{z}\leq\pi/4\}$,
\item $R_{3}$ be the extension of $\frac{-\bar{r}}{1-|r|^{2}}$ in the sector $\Omega_{3}:\{z:3\pi/4<\arg{z}\leq\pi\}$,
\item $R_{4}$ be the extension of $\frac{r}{1-|r|^{2}}$ in the sector $\Omega_{4}:\{z:\pi<\arg{z}\leq 5\pi/4\}$,
\item $R_{6}$ be the extension of $\bar{r}$ in the sector $\Omega_{6}:\{z:7\pi/4\leq\arg{z}< 2\pi\}$.
\end{enumerate}

Next we open sectors and define $A$ in the regions as follows
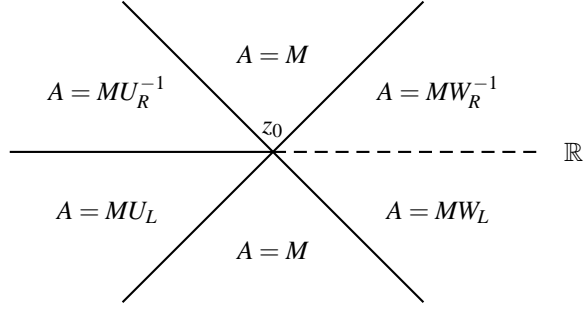
\begin{figure}[ht]
\begin{center}
\begin{pspicture}(-2,-2)(2,2)
\rput(0,0){\psline(-2,-2)(2,2)}
\rput(0,0){\psline(-2,2)(2,-2)}
\rput(0,0){\psline(0,0)(-3.5,0)}
\psset{linestyle=dashed}
\rput(0,0){\psline(0,0)(3.5,0)}
\rput(0,0.3){$z_{0}$}
\rput(4,0){$\mathbb R$}
\rput(0,1.3){$A=M$}
\rput(2.2,0.8){$A=MW_{R}^{-1}$}
\rput(2.2,-0.8){$A=MW_{L}$}
\rput(0,-1.3){$A=M$}
\rput(-2.2,-0.8){$A=MU_{L}$}
\rput(-2.2,0.8){$A=MU_{R}^{-1}$}
\end{pspicture}
\caption{The definition of $A$ in the different sectors $\Omega_{i}, i=1, \ldots, 6$.}
\end{center}
\end{figure}
where we use the extensions $R_{i}$ defined above for $U_{R}$, $U_{L}$, $W_{R}$ and $W_{L}$ in the sectors. The matrix $A$ now has the diagonal jump $U_{0}=(1-|r|^{2})^{\sigma_{3}}$ on $(-\infty,z_{0}]$ where
\begin{equation}
\sigma_{3}=\left(
  \begin{array}{cc}
    1 & 0 \\
0 &  -1
  \end{array}
\right).
\end{equation}
We can remove the jump across $(-\infty,z_{0}]$ by defining
\begin{equation}
B=A\delta^{-\sigma_{3}}.
\end{equation}
It is easily checked algebraically that $B$ so-defined only has jumps across the four diagonal rays pictured below.
\bigskip
\begin{figure}[ht]
\begin{center}
\begin{pspicture}(-2,-2)(2,2)
\rput(0,0){\psline(-2,-2)(2,2)}
\rput(0,0){\psline(-2,2)(2,-2)}
\rput(0,0.3){$z_{0}$}
\rput(3,1){
$\left(
  \begin{array}{cc}
    1 & 0 \\
R_{1}e^{2it\theta}\delta^{-2} & 1
  \end{array}
\right)$
}
\rput(3,-1){
$\left(
  \begin{array}{cc}
    1 & -R_{6}e^{-2it\theta}\delta^{2} \\
0 & 1
  \end{array}
\right)$
}
\rput(-3,-1){
$\left(
  \begin{array}{cc}
    1 & 0 \\
R_{4}e^{2it\theta}\delta^{-2} & 1
  \end{array}
\right)$
}
\rput(-3,1){
$\left(
  \begin{array}{cc}
    1 & -R_{3}e^{-2it\theta}\delta^{2} \\
0 & 1
  \end{array}
\right)$
}
\end{pspicture}
\caption{The jump matrices for $B$.}
\end{center}
\end{figure}
\subsubsection{Specify extensions to arrive at a pure $\dbar$--problem}

It is now time to return to the extensions $R_{i}$ and specify them explicitly. We seek $R_{1}$ satisfying
\begin{equation}
\begin{cases}
R_{1}=r \quad \textrm{on the ray $\{ z: \ z = r ,\ r \ge 0 \}$}\\
R_{1}=f_{1}:=\hat{r}_{0}(z-z_{0})^{-2i\nu}\delta^{2} \quad \textrm{on the ray $\{ z: \ z = r e^{i \pi / 4},\ r \ge 0 \}$}\\
|\dbar R_{1}|\leq c_{1}|z-z_{0}|^{-1/2}+c_{2}|r'| \quad \textrm{in $\{ z: \ z = r e^{i \theta},\ r \geq 0,\ 0\leq\theta\leq\pi/4 \}$}
\label{r1cond},
\end{cases}
\end{equation}

$R_{3}$ satisfying
\begin{equation}
\begin{cases}
R_{3}=\frac{-\bar{r}}{1-|r|^{2}} \quad \textrm{on the ray $\{ z: \ z = -r ,\ r \ge 0 \}$}\\
R_{3}=f_{3}:=\frac{-\bar{\hat{r}}_{0}}{1-|\hat{r}_{0}|^{2}}(z-z_{0})^{2i\nu}\delta^{-2} \quad \textrm{on the ray $\{ z: \ z = r e^{3i \pi / 4},\ r \ge 0 \}$}\\
|\dbar R_{3}|\leq c_{1}|z-z_{0}|^{-1/2}+c_{2}|r'|  \quad \textrm{in $\{ z: \ z = r e^{i \theta},\ r \geq 0,\ 3\pi/4\leq\theta\leq\pi \}$}
\label{r3cond},
\end{cases}
\end{equation}

$R_{4}$ satisfying
\begin{equation}
\begin{cases}
R_{4}=\frac{r}{1-|r|^{2}} \quad \textrm{on the ray $\{ z: \ z = r e^{5i \pi / 4},\ r \ge 0 \}$}\\
R_{4}=f_{4}:=\frac{\hat{r}_{0}}{1-|\hat{r}_{0}|^{2}}(z-z_{0})^{-2i\nu}\delta^{2} \quad \textrm{on the ray $\{ z: \ z = -r,\ r \ge 0 \}$}\\
|\dbar R_{4}|\leq c_{1}|z-z_{0}|^{-1/2}+c_{2}|r'| \quad \textrm{in $\{ z: \ z = r e^{i \theta},\ r \geq 0,\ \pi\leq\theta\leq 5\pi/4 \}$}
\label{r4cond},
\end{cases}
\end{equation}

and $R_{6}$ satisfying
\begin{equation}
\begin{cases}
R_{6}=\bar{r} \quad \textrm{on the ray $\{ z: \ z = r ,\ r \ge 0 \}$}\\
R_{6}=f_{6}:=\bar{\hat{r}}_{0}(z-z_{0})^{2i\nu}\delta^{-2} \quad \textrm{on the ray $\{ z: \ z = r e^{7i \pi / 4},\ r \ge 0 \}$}\\
|\dbar R_{6}|\leq c_{1}|z-z_{0}|^{-1/2}+c_{2}|r'| \quad \textrm{in $\{ z: \ z = r e^{i \theta},\ r \geq 0,\ 7\pi/4\leq\theta\leq 2\pi \}$}
\label{r6cond},
\end{cases}
\end{equation}
The reason for requiring the above conditions will be apparent soon.
\begin{prop}
There exist $R_{i}, i=1, 3, 4, 6$ satisfying conditions~(\ref{r1cond})-(\ref{r6cond}).
\end{prop}
\begin{proof}
For conciseness we give the proof for $R_{1}$ only; the corresponding arguments for the other $R_{i}$ follow immediately. Define
\[
R_{1}(u,v)=b(\arg(u+iv))r(u)+(1-b(\arg(u+iv)))f_{1}(u+iv)
\]
where
\[
b(\theta)=\cos(2\theta).
\]
It follows that $R_{1}$ satisfies the first two conditions in~(\ref{r1cond}) immediately. Note that
\[
\dbar R_{1}=(r-f_{1})\dbar b+\frac{b}{2}r',
\]
so
\[|\dbar R_{1}|\leq \frac{c_{1}}{|z-z_{0}|}\left[|r-r(z_{0})|+|r(z_{0})-f_{1}|\right]+c_{2}|r'|.\]
Note also that
\[|r-r(z_{0})|=\left\vert\int_{z_{0}}^{z} r'd\,s\right\vert \leq \int_{z_{0}}^{z}|r'|\,|d\,s| \leq ||r||_{L^{2}((z,z_{0}))}\cdot||1||_{L^{2}((z,z_{0}))}\leq c|z-z_{0}|^{1/2}.\]
Using the expansion in (\ref{gammaexpansion}), we also obtain that
\begin{eqnarray*}
f_{1}&=&\hat{r}_{0}(z-z_{0})^{-2i\nu}\delta^{2}\\
&=& \hat{r}_{0}\exp[2i\nu+2\beta(z_{0},z_{0})]\exp[2(\beta(z,z_{0})-\beta(z_{0},z_{0}))]\times\\
&\qquad&\exp[2i\nu((z-z_{0})\ln(z-z_{0})-(z-z_{0}+1)\ln(z-z_{0}+1))].
\end{eqnarray*}
We choose the free parameter $\hat{r}_{0}$ in such a way that
\begin{equation}\hat{r}_{0}\exp[2i\nu+2\beta(z_{0},z_{0})]=r(z_{0}).\label{rhat0}\end{equation}
Thus
\begin{equation*}
r(z_{0})-f_{1}= r(z_{0})-r(z_{0})(\exp[2i\nu((z-z_{0})\ln(z-z_{0})-(z-z_{0}+1)\ln(z-z_{0}+1))+2(\beta(z,z_{0})-\beta(z_{0},z_{0}))].
\end{equation*}
From (\ref{betaestimate1}) and (\ref{betaestimate2}) we obtain, uniformly in $\left\{z:z=re^{i\theta}, r>0, 0<\theta<\pi/4\right\}$, that
\[|\beta(z,z_{0})-\beta(z_{0},z_{0})|=\mathcal{O}(\sqrt{|z-z_{0}|}),\]
 and
\[|(z-z_{0})\ln(z-z_{0})|\leq\mathcal{O}(\sqrt{|z-z_{0}|}).\]
Therefore
\begin{eqnarray*}
|r(z_{0})-f_{1}|&=& r(z_{0})\left\{1-\exp[\mathcal{O}(\sqrt{|z-z_{0}|})]\right\}\\
&=& r(z_{0})\left\{\mathcal{O}(\sqrt{|z-z_{0}|})\right\}.
\end{eqnarray*}
Combining these estimates yields
\begin{equation}
|\dbar R_{1}| \leq c_{1}|z-z_{0}|^{-1/2} + c_{2}|r'|.
\label{dbarR1estimate}\end{equation}
\end{proof}

Recall that in the definition of the parabolic cylinder RHP (\ref{modrhp2}) there is a free parameter $r_{0}$. We set
\begin{equation}r_{0}=\hat{r}_{0}e^{i\nu\ln(8t)-4itz_{0}^{2}}\label{r0}.\end{equation}
Then by construction (conditions~(\ref{r1cond})-(\ref{r6cond})), the jumps for $B$ exactly match the the jumps for the parabolic cylinder RHP (\ref{modrhp2}). Thus we can define $E=BP^{-1}(\sqrt{8t}(z-z_{0}))$ which will have no jumps in the plane, but is a solution of the following problem.

\begin{dbp}
Find a 2$\times 2$ matrix--valued function $E$ satifying:
\begin{itemize}
\item $E$ is continuous in $\mathbb C$,
\item $E\to I$ as $z\to\infty$,
\item $\dbar E=E W$,
\end{itemize}
with
\[
W=
\begin{cases}
P\pmtwo{0}{0}{\dbar R_{1}e^{2it\theta}\delta^{-2}}{0}P^{-1}&\qquad\textrm{in sector $\Omega_{1}$},\\
P\pmtwo{0}{-\bar{\partial}R_{3}e^{-2it\theta}\delta^{2}}{0}{0}P^{-1}&\qquad\textrm{in sector $\Omega_{3}$},\\
P\pmtwo{0}{0}{\bar{\partial}R_{4}e^{2it\theta}\delta^{-2}}{0}P^{-1}&\qquad\textrm{in sector $\Omega_{4}$},\\
P\pmtwo{0}{-\bar{\partial}R_{6}e^{-2it\theta}\delta^{2}}{0}{0}P^{-1}&\qquad\textrm{in sector $\Omega_{6}$}\\
\pmtwo{0}{0}{0}{0} &\qquad\textrm{otherwise.}
\end{cases}
\]
\label{dbarprob}
\end{dbp}
The above problem is equivalent to the following integral equation.
\begin{equation}E=I-\frac{1}{\pi}\int\int\frac{EW}{s-z}\,d\,A(s)\label{integraleq}\end{equation}

\subsection{Analysis of the integral equation}\label{sect:estimates}
The integral equation~(\ref{integraleq}) can be written in operator notation as follows.
\begin{equation}\left[ \large{\textbf{1}} -J\right](E)=I,\label{eneumann}\end{equation}
where \begin{equation}J(E)=\frac{1}{\pi}\iint\frac{EW}{s-z}\,d\,A(s).\end{equation}
In order to invert the integral equation, we need to show that $J$ is small in norm. The following proposition gives the needed estimate.
\begin{prop}
\par There is a constant $c>0$ such that for all $t>0$, the following estimate holds.
\begin{equation}||J||_{L^{\infty}\to L^{\infty}}\leq c\,t^{-1/4}\label{jestimate}\end{equation}
\end{prop}
\begin{proof}
We give the details for sector $\Omega_{1}$ only as the corresponding arguments for the other sectors are identical with appropriate modifications.
Let $H\in L^{\infty}(\Omega_{1})$. Then
\begin{eqnarray*}
  |J(H)| &\leq &  \iint_{\Omega_{1}}\frac{|H\dbar R_{1}\delta^{-2}e^{2it\theta}|}{|s-z|}\,d\,A(s) \\
   &\leq & ||H||_{L^{\infty}(\omega_{1})}||\delta^{-2}||_{L^{\infty}(\omega_{1})} \iint_{\Omega_{1}}\frac{|\dbar R_{1}e^{2it\theta}|}{|s-z|}\,d\,A(s).
\end{eqnarray*}
Thus (\ref{dbarR1estimate}) yields
\[
  |J(H)| \leq C(I_{1}+I_{2}),
\]
where
\[I_{1} = \iint_{\Omega_{1}}\frac{|r'|e^{-tuv}}{|s-z|}\,d\,A(s),\]
and
\[I_{2} =  \iint_{\Omega_{1}}\frac{|s-z_{0}|^{-1/2}e^{-tuv}}{|s-z|}\,d\,A(s).\]
Since $r\in H^{1,1}(\mathbb R)$ we have
\begin{eqnarray*}
  |I_{1}| &\leq& \int_{0}^{\infty}\int_{v}^{\infty}\frac{|r'|e^{-tuv}}{|s-z|}\,d\,u\,d\,v  \\ &\leq&  \int_{0}^{\infty} e^{-tv^{2}}\int_{v}^{\infty}\frac{|r'|}{|s-z|}\,d\,u\,d\,v
   \leq  C\int_{0}^{\infty} e^{-tv^{2}}\bigg{\vert}\bigg{\vert}\frac{1}{s-z}\bigg{\vert}\bigg{\vert}_{L^{2}((v,\infty))}d\,v .\\
\end{eqnarray*}
Moreover
\begin{eqnarray*}
  \bigg{\vert}\bigg{\vert}\frac{1}{s-z}\bigg{\vert}\bigg{\vert}_{L^{2}((v,\infty))} &\leq & \left(\int_{\mathbb R}\frac{1}{|s-z|^{2}}\,du\right)^{1/2} \leq  \left(\frac{\pi}{|v-\beta|}\right)^{1/2}.
\end{eqnarray*}
Thus
\[
|I_{1}| \leq  C\int_{0}^{\infty} \frac{e^{-tv^{2}}}{\sqrt{\beta-v}}d\,v  \leq  C\left[\int_{0}^{\beta} \frac{e^{-tv^{2}}}{\sqrt{\beta-v}}d\,v+\int_{\beta}^{\infty} \frac{e^{-tv^{2}}}{\sqrt{v-\beta}}d\,v. \right]
\]
Using the fact $\sqrt{\beta}e^{-t\beta^{2}w^{2}}\leq ct^{-1/4}w^{-1/2}$, we obtain
\[
\int_{0}^{\beta} \frac{e^{-tv^{2}}}{\sqrt{\beta-v}}d\,v\leq \int_{0}^{1} \sqrt{\beta}\frac{e^{-t\beta^{2}w^{2}}}{\sqrt{1-w}}d\,w \leq \frac{c}{t^{1/4}}\int_{0}^{1} \frac{1}{\sqrt{w(1-w)}}d\,w \leq\frac{c}{t^{1/4}},
\]
whereas
\[
\int_{\beta}^{\infty} \frac{e^{-tv^{2}}}{\sqrt{v-\beta}}d\,v\leq \int_{0}^{\infty} \frac{e^{-tw^{2}}}{\sqrt{w}}d\,w\leq \frac{1}{t^{1/4}}\int_{0}^{\infty}    \frac{e^{-\lambda^{2}}}{\sqrt{\lambda}}d\,v\leq \frac{c}{t^{1/4}}.
\]
Hence
\[|I_{1}|\leq \frac{c}{t^{1/4}}.\]
To arrive at a similar estimate for $I_{2}$, we start with the following $L^{p}$--estimate for $p>2$.
\begin{eqnarray*}
\bigg{\vert}\bigg{\vert}\frac{1}{\sqrt{|s-z_{0}|}}\bigg{\vert}\bigg{\vert}_{L^{p}((d\,u))} &=& \left(\int_{z_{0}+v}^{\infty}\frac{1}{|u+iv-z_{0}|^{p/2}}\,d\,u\right)^{1/p} \\
&=& \left(\int_{v}^{\infty}\frac{1}{|u+iv|^{p/2}}\,d\,u\right)^{1/p} \\
&=& \left(\int_{v}^{\infty}\frac{1}{(u^{2}+v^{2})^{p/4}}\,d\,u\right)^{1/p} \\
&=& v^{(1/p-1/2)}\left(\int_{1}^{\infty}\frac{1}{(1+x^{2})^{p/4}}\,d\,x\right)^{1/p} \\
& \leq & c\,v^{1/p-1/2}
\end{eqnarray*}
Similarly to the $L^{2}$--estimate above, we obtain for $L^{q}$ where $1/q+1/p=1$.
\begin{eqnarray*}
  \bigg{\vert}\bigg{\vert}\frac{1}{s-z}\bigg{\vert}\bigg{\vert}_{L^{2}((v,\infty))} &\leq & C\,|v-\beta|^{1/q-1}
\end{eqnarray*}
It follows that
\[
|I_{2}| \leq  C \left[\int_{0}^{\beta} e^{-tv^{2}}v^{1/p-1/2}|v-\beta|^{1/q-1}d\,v + \int_{\beta}^{\infty} e^{-tv^{2}}v^{1/p-1/2}|v-\beta|^{1/q-1}d\,v \right].
\]
The first integral is handled similarly as above to yield the estimate
\[\int_{\beta}^{0} e^{-tv^{2}}v^{1/p-1/2}|v-\beta|^{1/q-1}d\,v\leq c\,t^{-1/4}.\]
To handle the second integral, we let $v=\beta+w$ so that the integral becomes
\[
\int_{0}^{\infty} e^{-t(\beta+w)^{2}}(\beta+w)^{1/p-1/2}w^{1/q-1}d\,w\leq \int_{0}^{\infty} e^{-tw^{2}}(w)^{1/p-1/2}w^{1/q-1}d\,w.
\]
Using the substitution $y=\sqrt{t}w$ then yields
\[\int_{\beta}^{\infty} e^{-tv^{2}}v^{1/p-1/2}|v-\beta|^{1/q-1}d\,v\leq c\,t^{-1/4}.\]
Combining the previous estimates we obtain
\[|I_{2}| \leq  c\,t^{-1/4},\]
and the result follows.
\end{proof}
For $t$ sufficiently large the integral equation (\ref{eneumann}) may be inverted by Neumann series. It follows that
\begin{equation}E=I+\mathcal{O}\left(t^{-1/4}\right).\label{E}\end{equation}
\subsection{Asymptotics for NLS}
Recall that
\[M=EP(\sqrt{8t}(z-z_{0}))\delta^{\sigma_{3}}=\left(I+\frac{E_{1}}{z}+\mathcal{O}\left({z^{-2}}\right)\right)\left(I+\frac{P_{1}^{\infty}}{\sqrt{8t}(z-z_{0})}+\mathcal{O}\left({(\sqrt{8t}(z-z_{0}))^{-2}}\right)\right)\left(I+\ldots\right)\]
for $z$ in the sectors $\Omega_{2}$ and $\Omega_{5}$. Denote by $M_{1}$ the coefficient of $z^{-1}$ in the Laurent expansion of $M$
\[M=I+\frac{M_{1}}{z}+\ldots,\]
so that
\[M_{1}=E_{1}+\frac{P_{1}^{\infty}}{\sqrt{8t}}+\pmtwo{\delta_{1}}{0}{0}{-\delta_{1}}.\]
Hence the contribution to the off-diagonal entries of $M_{1}$ which are relevant for the asymptotics (recall \ref{qMrelation}) come from $P_{1}^{\infty}$ and $E_{1}$.  From the integral equation~(\ref{integraleq}) satisfied by $E$, we have
\begin{equation}E_{1}=\frac{1}{\pi}\iint EW\,d\,A,\label{E1}\end{equation}
and one may verify, using estimates of the type employed in section~\ref{sect:estimates}, that
\[|E_{1}|\leq c\,t^{-3/4}.\]
Specifically
\[|E_{1}|\leq\int_{0}^{\infty}\int_{v}^{\infty}e^{-tuv}\,|r'|\,d\,u\,d\,v+\int_{0}^{\infty}\int_{z_{0}+v}^{\infty}e^{-tuv}\,|z-z_{0}|^{1/2}\,d\,u\,d\,v=I_{3}+I_{4},\]
and the Cauchy-Schwarz inequality gives
\[I_{3}\leq C \int_{0}^{\infty}\left(\int_{v}^{\infty}e^{-tuv}\,d\,u\right)^{1/2}\,d\,v\leq C\,t^{-3/4}.\]
Similarly, the H\"older inequality for $1/p+1/q=1$,  $2<p<4$ yields
\[I_{4}\leq \int_{0}^{\infty}v^{1/p-1/2}\left(\int_{v}^{\infty}e^{-qtuv}\,d\,u\right)^{1/q}\,d\,v\leq\frac{c}{t^{1/q}}\int_{0}^{\infty}v^{2/p-3/2}e^{-tv^{2}}\,d\,v\leq C\,t^{-3/4}.\]
The main result Theorem~(\ref{mainresult}) then follows from
\begin{equation}q(x,t) =2i\frac{\left(P_{1}^{\infty}\right)_{12}}{\sqrt{8t}}+\mathcal{O}\left(t^{-3/4}\right),\end{equation}
where $\left(P_{1}^{\infty}\right)_{12}$ is given in (\ref{p1infinty}) and written explicitly in terms of $r(z_{0})$ using (\ref{r0}), (\ref{rhat0}) and (\ref{beta}). The identity \[\mid\Gamma(i\nu)\mid^{2}=\mid\Gamma(-i\nu)\mid^{2}=\pi\left(\sinh(\pi\nu)\right)^{-1}\] is useful in carrying out the algebraic manipulations needed.
\begin{remark}
It is interesting to observe that the calculations we have presented actually yield an asymptotic expansion (of sorts) for $q(x,t)$. For example, instead of the integral estimates that follow (\ref{E1}), one could use (\ref{E}) in (\ref{E1}), yielding
\begin{equation}E_{1}=\frac{1}{\pi}\iint W\,d\,A + \mathcal{O}\left(t^{-1}\right),\end{equation}
which in turn implies
\begin{equation}q(x,t) =2i\frac{\left(P_{1}^{\infty}\right)_{12}}{\sqrt{8t}} +\frac{2i}{\pi}\left(\iint W\,d\,A\right)_{12}+ \mathcal{O}\left(t^{-1}\right).\end{equation}
More generally, one may use the Neumann expansion for $E$ arising from (\ref{eneumann}) to obtain an asymptotic expansion to arbitrary order. Extraction of more detail from the subsequent terms in the asymptotic expansion involves an analysis of explicit double integrals built out of the reflection coefficient and the parabolic cylinder functions.
\end{remark}


\end{document}